\documentclass{article}
\usepackage{latexsym}
\usepackage{amstext}
\usepackage{amsmath}
\usepackage{amssymb}
\usepackage{amsthm}

\newtheorem{theorem}{Theorem}
\newtheorem{theoremb}{Theorem}

\newtheorem{corollarya}[theorem]{Corollary}
\newtheorem{corollaryb}[theoremb]{Corollary}

\newtheorem{lemmaa}[theorem]{Lemma}
\newtheorem{lemmab}[theoremb]{Lemma}

\newtheorem{theorema}[theorem]{Theorem}

\theoremstyle{remark}

\renewcommand{\thetheorem}{\arabic{theorem}A}

\def \co {\mathcal{O}}
\def \ch {\mathcal{H}}
\def \cl {\mathcal{L}}
\def \cm {\mathcal{M}}
\def \cpn {\mathbb{P}^n_{\mathbb{C}}}
\def \kpn {\mathbb{P}^n_k}
\def \Kpn {\mathbb{P}^n_K}
\DeclareMathOperator{\Span}{Span}

\begin{document}
\bibliographystyle{amsplain}
\title{The Dimensions of Integral Points and Holomorphic Curves on the Complements of Hyperplanes}
\author{Aaron Levin}
\date{}          
\maketitle
\begin{abstract}
In this article we completely determine the possible dimensions of integral points and holomorphic curves on the complement of a union of hyperplanes in projective space.  Our main theorems generalize a result of Evertse and Gy\"ory, who determined when all sets of integral points (over all number fields) on the complement of a union of hyperplanes are finite, and a result of Ru, who determined when all holomorphic maps to the complement of a union of hyperplanes are constant.  The main tools used are the $S$-unit lemma and its analytic analogue, Borel's lemma.
\end{abstract}
\section{Introduction}
In this article we completely determine the possible dimensions of a set of integral points or a holomorphic curve on the complement of a union of hyperplanes in projective space.  When we say determine, we mean that from our results one may construct an algorithm that takes as input the defining equations for the hyperplanes and returns as output, after a finite number of steps, the possible dimensions of a set of integral points or a holomorphic curve on the complement of the hyperplanes. In Vojta's Nevanlinna-Diophantine dictionary, a holomorphic curve in a complex variety $V$, i.e., a holomorphic map $f:\mathbb{C}\to V$, corresponds to a set of integral points on $V$, assuming $V$ is defined over a number field.  We therefore expect that the possible dimensions of a holomorphic curve and a set of integral points on $V$ will be the same, and indeed, we will see that this is the case when $V$ is the complement of a union of hyperplanes. 

Our main theorems generalize a result of Evertse and Gy\"ory \cite{Ev}, who determined when all sets of integral points (over all number fields) on the complement of a union of hyperplanes are finite, and a result of Ru \cite{Ru2}, who determined when all holomorphic maps to the complement of a union of hyperplanes are constant.  We also generalize results of Fujimoto \cite{Fu2} and Green \cite{Gr}, who bounded the dimension of a holomorphic curve on the complement of hyperplanes in general position.
\section{Results}
If $\ch$ is a set of hyperplanes on $\mathbb{P}^n$, we denote the union of the hyperplanes by $|\ch|$.  We denote by $\cl$ a set of linear forms in $x_0,\ldots, x_n$ which define the hyperplanes in $\ch$.  We let $(\cl)$ be the complex vector space generated by the elements of $\cl$.

Let $k$ be a number field and $M_k$ a complete set of inequivalent places of $k$.  Let $S$ be a finite set of places in $M_k$ (which we will always assume contains the archimedean places $S_\infty$ of $k$).  Let $\co_{k,S}$ denote the ring of $S$-integers of $k$.  If the set $|\ch|$ is defined over $k$, we call a set $R\subset \kpn\backslash |\ch|$ a set of $S$-integral points on $\kpn \backslash |\ch|$ if for every regular function $f$ on $\kpn \backslash |\ch|$ there exists $a\in k^*$ such that $af(P)\in \co_{k,S}$ for all $P\in R$.  Equivalently, $R$ is a set of $S$-integral points on $\kpn \backslash |\ch|$ if there exists an affine embedding $\kpn \backslash |\ch|\subset \mathbb{A}^N_k$ such that every $P\in R$ has $S$-integral coordinates.

We now state our two main theorems.
\begin{theorema}
\label{maina}
Let $\ch$ be a set of hyperplanes in $\kpn$.  Let $\cl$ be a corresponding set of linear forms.  Let $m=\dim \cap_{H\in   \ch}H$.  Then there exists a number field $K\supset k$, $S\subset M_K$, and a set of $S$-integral points $R$ on $\Kpn\backslash |\ch|$ with $\dim R=m+1$.  There exists a number field $K\supset k$, $S\subset M_K$, and a set of $S$-integral points $R$ on $\Kpn \backslash |\ch|$ with $\dim R=d>m+1$ if and only if there exists a partition of $\cl$ into $d-m$ nonempty disjoint subsets $\cl_i$,
\begin{equation}
\label{h1}
\cl=\bigsqcup_{i=1}^{d-m} \cl_i, \quad \cl_i\neq \emptyset \text{ for all }i,
\end{equation}
such that
\begin{equation}
\label{h2}
\cl\cap \sum_{j=1}^{d-m}\left((\cl_j)\cap \sum_{i\neq j}(\cl_i)\right)=\emptyset.
\end{equation}
\end{theorema}
\begin{theoremb}
\label{mainb}
Let $\ch$ be a set of hyperplanes in $\cpn$.  Let $\cl$ be a corresponding set of linear forms.  Let $m=\dim \cap_{H\in   \ch}H$.  Then there exists a holomorphic map $f:\mathbb{C}\to \cpn \backslash |\ch|$ with $\dim f(\mathbb{C})=m+1$.  There exists a holomorphic map $f:\mathbb{C}\to \cpn \backslash |\ch|$ with $\dim f(\mathbb{C})=d>m+1$ if and only if there exists a partition of $\cl$ into $d-m$ nonempty disjoint subsets $\cl_i$ satisfying \eqref{h1} and \eqref{h2} above.
\end{theoremb}
Note that in the these theorems and elsewhere we define $\dim \emptyset=-1$.
Given a partition of $\cl$, one needs only elementary linear algebra to check (\ref{h2}).  So running over all partitions of $\cl$ and checking (\ref{h2}), we may determine the possible dimensions of a set of $S$-integral points on $\Kpn \backslash |\ch|$ (over all $K$ and $S$) or of a holomorphic map $f:\mathbb{C} \to \cpn\backslash |\ch|$.  

We will prove Theorems \ref{maina} and \ref{mainb} in the next two sections.  We now mention some consequences of these theorems.  The next two corollaries are immediate from the theorems.
\begin{corollarya}[Evertse and Gy\"ory \cite{Ev}]
Let $\ch$ be a set of hyperplanes in $\kpn$.  Let $\cl$ be a corresponding set of linear forms.  Then all sets of $S$-integral points on $\Kpn\backslash |\ch|$ are finite for every number field $K\supset k$ and $S\subset M_K$ if and only if $\cap_{H\in \ch}H=\emptyset$ and for every proper nonempty subset $\cl_1$ of $\cl$, 
\begin{equation*}
\cl \cap(\cl_1)\cap(\cl\backslash \cl_1)\neq \emptyset.
\end{equation*}

\end{corollarya}
\begin{corollaryb}[Ru \cite{Ru2}]
Let $\ch$ be a set of hyperplanes in $\cpn$.  Let $\cl$ be a corresponding set of linear forms.  Then all holomorphic maps $f:\mathbb{C}\to \cpn \backslash |\ch|$ are constant (i.e., $\cpn \backslash |\ch|$ is Brody hyperbolic) if and only if $\cap_{H\in \ch}H=\emptyset$ and for every proper nonempty subset $\cl_1$ of $\cl$, 
\begin{equation*}
\cl \cap(\cl_1)\cap(\cl\backslash \cl_1)\neq \emptyset.
\end{equation*}
\end{corollaryb}
Let $[x]$ denote the greatest integer in $x$.
\begin{corollarya}
Let $\ch$ be a set of hyperplanes in $\kpn$.  Suppose that the intersection of any $s+1$ distinct hyperplanes in $\ch$ is empty.  Let $r=\#\ch$.  Suppose $r>s$.  Then for every number field $K\supset k$ and $S\subset M_K$, for all sets $R$ of $S$-integral points on $\Kpn\backslash |\ch|$,
\begin{equation}
\label{ineqa}
\dim R \leq \left[\frac{s}{r-s}\right].
\end{equation}
In particular, if $r>2s$, then all such $R$ are finite.  Furthermore, if the hyperplanes in $\ch$ are in general position ($s=n$), then the bound in \eqref{ineqa} is achieved by some $R$.
\end{corollarya}
\begin{corollaryb}
\label{Green}
Let $\ch$ be a set of hyperplanes in $\cpn$.  Suppose that the intersection of any $s+1$ distinct hyperplanes in $\ch$ is empty.  Let $r=\#\ch$.  Suppose $r>s$.  Then for all holomorphic maps $f:\mathbb{C}\to \cpn \backslash |\ch|$, 
\begin{equation}
\label{ineqb}
\dim f(\mathbb{C})\leq \left[\frac{s}{r-s}\right].
\end{equation}
In particular, if $r>2s$ then $\cpn\backslash |\ch|$ is Brody hyperbolic.  Furthermore, if the hyperplanes in $\ch$ are in general position ($s=n$), then the bound in \eqref{ineqb} is achieved by some $f$.
\end{corollaryb}
Corollary \ref{Green} generalizes theorems of Fujimoto \cite{Fu2} and Green \cite{Gr}, who each independently proved the case when the hyperplanes are in general position.  Working in a different direction, Noguchi and Winkelmann \cite{No} have generalized Fujimoto and Green's result (and its arithmetic analogue) to hypersurfaces of projective space in  general position.

We now prove both corollaries simultaneously.
\begin{proof}
To prove (\ref{ineqa}) and (\ref{ineqb}), let $d>\frac{s}{r-s}$ be an integer.  Let $\cl=\sqcup_{i=1}^{d+1}\cl_i$ be a partition of $\cl$ into nonempty disjoint subsets $\cl_i$, $i=1,\ldots, d+1$.  Note that since $r>s$, $\cap_{H\in \ch}H=\emptyset$, so $m=-1$ in Theorems \ref{maina} and \ref{mainb}.  Using $d>\frac{s}{r-s}$, we see that there exists an index $i_0$ such that $\#\cup_{i\neq i_0}\cl_i\geq \frac{d}{d+1}r>s$.  Since at most $s$ of the hyperplanes meet at a point, we therefore have that $\sum_{i\neq i_0}(\cl_i)$ is the whole $n+1$-dimensional vector space of linear forms, and so in particular, we see that
\begin{equation*}
\cl\cap \sum_{j=1}^{d+1}\left((\cl_j)\cap \sum_{i\neq j}(\cl_i)\right)\neq\emptyset.
\end{equation*}
Thus, using Theorems \ref{maina} and \ref{mainb}, we obtain the desired inequalities.

For the last assertions of the corollaries, where $s=n$, let $d=\left[\frac{n}{r-n}\right]$.  Then we may partition $\cl$ into nonempty disjoint subsets $\cl_i$ for $i=1,\ldots, d+1$ such that $\#\cl_i \geq r-n$ for all $i$.  Suppose that (\ref{h2}) doesn't hold.  Let $L\in \cl$ with 
\begin{equation*}
L\in \sum_{j=1}^{d+1}\left((\cl_j)\cap \sum_{i\neq j}(\cl_i)\right).
\end{equation*}
For some index $i_0$, $L\in\cl_{i_0}$. Note that
\begin{equation*}
\sum_{j=1}^{d+1}\left((\cl_j)\cap \sum_{i\neq j}(\cl_i)\right)\subset \sum_{i\neq i_0}(\cl_i),
\end{equation*}
so that $L\in \sum_{i\neq i_0}(\cl_i)$.  From $\#\cl_{i_0}\geq r-n$, we have that $\#\cup_{i\neq i_0}\cl_i\leq n$.  This implies that the set of at most $n+1$ linear forms $\cup_{i \neq i_0}\cl_i\cup \{L\}$ is linearly dependent.  This contradicts the hypothesis that the hyperplanes were in general position.  Therefore (\ref{h2}) holds for this partition of $\cl$, and so we are done by Theorems \ref{maina} and \ref{mainb}.
\end{proof}

\section{A Reformulation of the Problem}
In this section we give a simple reformulation of the problem from the introduction.  If $Y\not\subset |\ch|$ is a linear subspace of $\mathbb{P}^n$ then we define 
\begin{equation*}
\ch|_Y=\{H\cap Y\mid H\in \ch\}.
\end{equation*}
Note that $\ch|_Y$ may contain fewer hyperplanes than $\ch$.  We will denote by $\cl|_Y$ a set of linear forms defining the hyperplanes in $\ch|_Y$.  

Consider the condition
\begin{equation}
\label{cond}
Y\subset \cpn \text{ is a linear space}, Y\not\subset |\ch|, \text{and } \cl|_Y \text{ is a linearly independent set}.
\end{equation}
We now reformulate our problem in terms of this condition.
\begin{theorema}
\label{refa}
Let $\ch$ be a set of hyperplanes in $\kpn$.  Let $\cl$ be a corresponding set of linear forms.  There exists a number field $K\supset k$, $S\subset M_K$, and a set of $S$-integral points $R$ on $\Kpn\backslash |\ch|$ with $\dim R=d$ if and only if there exists a $Y$ satisfying \eqref{cond} with $\dim Y=d$.
\end{theorema}
\begin{theoremb}
\label{refb}
Let $\ch$ be a set of hyperplanes in $\cpn$.  Let $\cl$ be a corresponding set of linear forms.  There exists a holomorphic map $f:\mathbb{C}\to \cpn \backslash |\ch|$ with $\dim f(\mathbb{C})=d$ if and only if there exists a $Y$ satisfying \eqref{cond} with $\dim Y=d$.
\end{theoremb}
We will see that these theorems are simple consequences of the following two fundamental lemmas.  We begin by giving the $S$-unit lemma.
\begin{lemmaa}[$S$-unit Lemma]
Let $k$ be a number field and let $n>1$ be an integer.  Let $\Gamma$ be a finitely generated subgroup of $k^*$.  Then all but finitely many solutions of the equation
\begin{equation}
u_0+u_1+\cdots+u_n=1, u_i\in \Gamma
\end{equation}
satisfy an equation of the form $\sum_{i\in I}u_i=0$, where $I$ is a subset of $\{0,\ldots,n\}$.
\end{lemmaa}
The analytic analogue of the $S$-unit lemma is Borel's lemma.
\begin{lemmab}[Borel's Lemma]
Let $f_1,\ldots,f_n$ be entire functions without zeroes on $\mathbb{C}$.  Suppose that
\begin{equation}
f_1+\cdots+f_n=1.
\end{equation}
Then $f_i$ is constant for some $i$.
\end{lemmab}

\begin{lemmaa}
\label{hypa}
Let $\ch$ be a set of hyperplanes in $\kpn$, and let $\cl$ be a corresponding set of linear forms.  There does not exist a Zariski-dense set of $S$-integral points on $\Kpn \backslash |\ch|$ for all number fields $K\supset k$ and $S\subset M_K$ if and only if $\cl$ is a linearly dependent set.  Furthermore, in this case any set $R$ of $S$-integral points on $\Kpn \backslash |\ch|$ is contained in a finite union of hyperplanes of $\Kpn$. 
\end{lemmaa}
\begin{proof}
Suppose $\cl$ is a linearly dependent set.  Let $\{L_1,\ldots,L_m\}\subset \cl$ be a minimal linearly dependent subset, that is, no proper subset is linearly dependent.  Then $\sum_{i=1}^{m-1} c_iL_i=c_mL_m$  for some choice of $c_i\in k^*$.  Let $R$ be a set of $S$-integral points on $\Kpn\backslash |\ch|$.
Since all of the poles of $\frac{L_i}{L_m}$ lie in $|\ch|$, there exists an $a\in K^*$ such that $af$ takes on $S$-integral values on $R$.  Since the poles of $\frac{L_m}{L_i}$ also lie in $|\ch|$, the same reasoning applies to $\frac{L_m}{L_i}$.  Therefore $\frac{L_i}{L_m}(P)$ lies in only finitely many cosets of $\co_{K,S}^*$ for $P\in R$.  
By enlarging $S$, we may assume without loss of generality that $\frac{c_iL_i}{c_mL_m}(P)$ is an $S$-unit for all $P\in R$ and $i=1,\ldots,m$.  Since $\sum_{i=1}^{m-1} \frac{c_iL_i}{c_mL_m}(P)=1$ for all $P\in R$, by the $S$-unit lemma, it follows that all $P\in R$ either belong to one of the hyperplanes defined by $\sum_{i\in I}c_iL_i=0$ for some subset $I\subset \{1,\ldots,m-1\}$ (this equation is nontrivial by the minimality of the linear dependence relation) or they belong to a hyperplane defined by $\frac{c_iL_i}{c_mL_m}=t\in T$, where $T \subset \co_{K,S}^*$ is a finite subset containing the elements that appear in the exceptional solutions to the $S$-unit equation $\sum_{i=1}^{m-1} x_i=1$.  Thus $R$ is contained in a finite union of hyperplanes of $\Kpn$.

Conversely, suppose that $\cl$ is a linearly independent set.  After a $k$-linear change of coordinates, we may assume that $\cl=\{x_0,\ldots,x_m\}$ for some $m\leq n$.  Let $K\supset k$ be a number field with $\co_K^*$ infinite.  Let $S$ be the set of archimedean places of $K$.  Let $R$ be the set of points in $\Kpn$ which have a representation where the coordinates are all ($S$-)units.  Then it is easy to see that $R$ is a set of $S$-integral points on $\Kpn\backslash |\ch|$.  We now show that $R$ is Zariski-dense in $\Kpn$.  This follows from the $S$-unit lemma, or the following more elementary argument.  Consider the set of homogeneous polynomials in $x_0,\ldots, x_n$ that vanish on $R\subset \Kpn$.  If this set is nonempty, let $p$ be a polynomial in this set with a minimal number of terms.  Let $x_i$ be a variable that appears with different powers in two monomials of $p$.  Let $u\in \co_K^*$ be a unit that is not a root of unity.  Let $q$ be the homogeneous polynomial obtained from $p$ by the substitution $x_i\mapsto ux_i$.  Then $q$ also vanishes on $R$.  By our choice of $u$ and $x_i$, $q$ is not a scalar multiple of $p$. However, $p$ and $q$ contain the same monomials.  Therefore, there exists a linear combination of $p$ and $q$ that vanishes on $R$ and has strictly fewer terms than $p$.  This contradicts the minimality of $p$, so $R$ is Zariski-dense in $\Kpn$.
\end{proof}
\begin{lemmab}
\label{hypb}
Let $\ch$ be a set of hyperplanes in $\cpn$, and let $\cl$ be a corresponding set of linear forms.  There does not exist a holomorphic map $f:\mathbb{C}\to \cpn \backslash |\ch|$ with Zariski-dense image if and only if $\cl$ is a linearly dependent set.  Furthermore, in this case all such holomorphic maps $f$ have image contained in a hyperplane of $\cpn$. 
\end{lemmab}
\begin{proof}
Suppose $\cl$ is a linearly dependent set.  Let $\{L_1,\ldots,L_m\}\subset\cl$ be a minimal linearly dependent set.  Then there exist nonzero constants $c_2,\ldots,c_m$ such that $\sum_{i=2}^m c_i \frac{L_i}{L_1}=1$.  Let $f:\mathbb{C}\to \cpn \backslash |\ch|$ be a holomorphic map and let $g_i=\frac{L_i}{L_1}\circ f$.  Then $g_i$ is an entire function without zeroes on $\mathbb{C}$ since the image of $f$ omits $|\ch|$. We also have $\sum_{i=2}^mc_ig_i=1$.  By Borel's lemma, $g_i=\alpha$ for some $i$ and some constant $\alpha \in \mathbb{C}$.  Therefore the image of $f$ is contained in the hyperplane defined by $L_i-\alpha L_1=0$.

Conversely, suppose that $\cl$ is a linearly independent set.  After a linear change of coordinates, we may assume that $\cl=\{x_0,\ldots,x_m\}$ for some $m\leq n$.  Let $f:\mathbb{C}\to \cpn \backslash |\ch|$ be defined by $f=(1,e^z,e^{z^2},\ldots,e^{z^n})$.  Looking at the growth as $z\to \infty$, it is clear that no homogeneous polynomial can vanish on $f(\mathbb{C})$, and so $f(\mathbb{C})$ is Zariski-dense in $\cpn$.
\end{proof}
Theorems \ref{refa} and \ref{refb} now follow rather directly from Lemmas \ref{hypa} and \ref{hypb}.
\begin{proof}[Proof of Theorems \ref{refa} and \ref{refb}]
We first make two general observations.  In condition (\ref{cond}), even in the arithmetic case, we have allowed complex linear spaces $Y$.  However, when the hyperplanes in $\ch$ are defined over a number field $k$, it is easily seen that
\begin{multline*}
\max\{\dim Y\mid Y \text{ satisfies }(\ref{cond})\} =\\
\max\{\dim Y\mid Y \text{ is defined over a number field and $Y$ satisfies }(\ref{cond})\}.
\end{multline*}

Secondly, if $Y$ satisfies (\ref{cond}), then for any $0\leq d'\leq \dim Y$ there exists a $Y'$ satisfying (\ref{cond}) with $\dim Y'=d'$.  To see this, let $\ch|_Y=\{H_1,\ldots,H_m\}$.  If $m\geq \dim Y+1-d'$, then let $Y'$ be any linear space $Y'\not\subset |\ch|$ with 
\begin{equation*}
\cap_{i=1}^{\dim Y+1-d'}H_i\subset Y'
\end{equation*}
and $\dim Y'=d'$.  Then $Y'$ satisfies (\ref{cond}).  If $m< \dim Y+1-d'$, then $\cap_{i=1}^{\dim Y+1-d'}H_i$ consists of a single linear space of dimension $\geq d'$.  Therefore, we may choose $Y'$ to be a linear space $Y'\not\subset |\ch|$ with $\dim Y'=d'$ and $\ch|_{Y'}$ consisting of a single hyperplane (ignoring the trivial case $d'=0$).

Suppose now that there exists a $Y$ satisfying (\ref{cond}). Using the remarks above, we see that if the hyperplanes in $\ch$ are defined over a number field, then there exists a $Y'$ satisfying (\ref{cond}) with $Y'$ defined over a number field and $\dim Y'=\dim Y$.  So, restricting things to $Y$ (or $Y'$), it is immediate from Lemmas \ref{hypa} and \ref{hypb} that there exists a number field $K\supset k$, $S\subset M_K$, and a set of $S$-integral points $R$ on $\Kpn \backslash |\ch|$ with $\dim R=\dim Y$, and there exists a holomorphic map $f:\mathbb{C}\to \cpn \backslash |\ch|$ with $\dim f(\mathbb{C})=\dim Y$.

Let $R$ be a set of $S$-integral points on $\Kpn \backslash |\ch|$ or the image of a holomorphic map $f:\mathbb{C}\to \cpn\backslash |\ch|$.  Repeatedly applying Lemma \ref{hypa} or \ref{hypb} (using that when $\cl$ is a linearly dependent set, $R$ is contained in a union of projective spaces to which the theorems may be applied again), we see that there exists a $Y$ satisfying (\ref{cond}) with $\dim R\leq \dim Y$.  By our earlier remarks, there then exists a $Y'$ satisfying (\ref{cond}) with $\dim Y'=\dim R$.
\end{proof}

\section{Proofs of Main Theorems}
We first make one more definition.  We define the zero set of a set of linear forms $\cl$ in $n+1$ variables to be the linear variety 
\begin{equation*}
Z(\cl)=\{P\in\cpn\mid L(P)=0 \text{ for all } L\in \cl\}.
\end{equation*}
Using Theorems \ref{refa} and \ref{refb} we are reduced to computing, for a given set of hyperplanes $\ch$, the possible dimensions of a linear space $Y$ satisfying (\ref{cond}).  Theorems \ref{maina} and \ref{mainb} are therefore equivalent to the following theorem.
\renewcommand{\thetheorem}{\arabic{theorem}}
\begin{theorem}
Let $\ch$ be a set of hyperplanes in $\cpn$.  Let $\cl$ be a corresponding set of linear forms.  Let $m=\dim \cap_{H\in   \ch}H$.  Then there exists a $Y$ satisfying \eqref{cond} with $\dim Y=m+1$.  There exists a $Y$ satisfying \eqref{cond} with $\dim Y=d>m+1$ if and only if there exists a partition of $\cl$ into $d-m$ nonempty disjoint subsets $\cl_i$ satisfying \eqref{h1} and \eqref{h2}.
\end{theorem}
\begin{proof}
We first prove our assertion about the existence of a $Y$ satisfying (\ref{cond}) with $\dim Y=m+1$.  If $m=-1$ this is trivial.  Otherwise, we may take $Y$ to be any linear subspace of $\cpn$ of dimension $m+1$ with $\cap_{H\in   \ch}H\subset Y\not \subset |\ch|$.  In this case $\cl|_Y$ consists of a single linear form, which is therefore a linearly independent set.

Suppose now that there exists a $Y$ satisfying (\ref{cond}) with $\dim Y=d>m+1$.  Let $m'=\dim \cap_{H'\in H|_Y}H'$.  Then $m'\leq m$.  Since $Y$ satisfies (\ref{cond}), $\ch|_Y$ consists of exactly $d-m'$ hyperplanes of $Y$, $H_1',\ldots,H_{d-m'}'$.  Let $\cl_i$ for $i=1,\ldots, d-m'$ consist of the linear forms in $\cl$ that define hyperplanes which intersect $Y$ in $H_i'$.  Then we get a partition of $\cl$ into $d-m'$ nonempty disjoint subsets $\cl_i$, $\cl=\sqcup_{i=1}^{d-m'} \cl_i$.  Let $j\in \{1,\ldots,d-m'\}$.  Then 
\begin{equation*}
Z\left(\sum_{i\neq j}(\cl_i)\right)\supset \bigcap_{i \neq j} H_i'.
\end{equation*}
It follows from the linear independence of $\cl|_Y$ that $\cap_{i \neq j} H_i'$ contains a point in $Y$ not contained in $H_j'$.  Therefore 
\begin{equation*}
Z\left((\cl_j)\cap \sum_{i\neq j}(\cl_i)\right)\supset \Span (H_j',\cap_{i \neq j} H_i')=Y.
\end{equation*}
So 
\begin{equation*}
Z\left(\sum_{j=1}^{d-m'}\left((\cl_j)\cap \sum_{i\neq j}(\cl_i)\right)\right)\supset Y.
\end{equation*}
Since $Y\not\subset |\ch|$, we must therefore have $\cl\cap \sum_{j=1}^{d-m'}\left((\cl_j)\cap \sum_{i\neq j}(\cl_i)\right)=\emptyset$.  Let $\cm_1=\cup_{i=1}^{m+1-m'}\cl_i$ and $\cm_i=\cl_{i+m-m'}$ for $i=2,\ldots,d-m$.  It is  straightforward to verify that 
\begin{equation*}
\sum_{j=1}^{d-m}\left((\cm_j)\cap \sum_{i\neq j}(\cm_i)\right)\subset\sum_{j=1}^{d-m'}\left((\cl_j)\cap \sum_{i\neq j}(\cl_i)\right).
\end{equation*}
Therefore (\ref{h1}) and (\ref{h2}) are satisfied (with $\cl_i=\cm_i$).

In the other direction, suppose that there exists a partition of $\cl$ into $d-m$ nonempty disjoint subsets $\cl_i$, $\cl=\sqcup_{i=1}^{d-m} \cl_i$, such that (\ref{h2}) is satisfied.  Let 
\begin{equation*}
U_0=\sum_{j=1}^{d-m}\left((\cl_j)\cap \sum_{i\neq j}(\cl_i)\right).
\end{equation*}
We now define vector spaces $U_i$ for $i=0,\ldots, d-m$ such that 
\begin{enumerate}
\item $U_i\subset U_j$ for $i<j$.
\item $\dim U_i\cap (\cl_i)=\dim (\cl_i)-1$ for $i>0$.
\item $U_i \cap \cl=\emptyset$.
\item $U_i=\sum_{j=1}^{d-m}U_i\cap (\cl_j)$.
\end{enumerate}
Clearly the space $U_0$ we have defined satisfies these properties.  

We inductively define the vector spaces $U_i$.  Suppose that we have defined a $U_{i-1}$ with the above properties.  Since $U_{i-1}\cap \cl=\emptyset$, it follows that $U_{i-1} \cap (\cl_i)$ is a proper subspace of $(\cl_i)$.  Therefore, since $\mathbb{C}$ is infinite, 
\begin{equation*}
\bigcup_{L\in \cl_i} (U_{i-1}\cap (\cl_i)+L)\neq (\cl_i).
\end{equation*}
So we easily see that we may add elements of $(\cl_i)$ to $U_{i-1}$ to get a space $U_i$ with $\dim U_i\cap (\cl_i)=\dim (\cl_i)-1$ and  $U_i \cap \cl_i=\emptyset$.  Also, we have $U_{i-1}\subset U_i$, and since we added only elements of $(\cl_i)$, we get from the corresponding property for $U_{i-1}$ that $U_i=\sum_{j=1}^{d-m}U_i\cap(\cl_j)$.  In fact, we have that 
\begin{equation}
\label{ueq}
U_i= U_i\cap(\cl_i)+ \sum_{j\neq i}U_{i-1}\cap(\cl_j).
\end{equation}
To show that $U_i$ satisfies all the required properties, it only remains to show that $U_i \cap \cl=\emptyset$.  

Suppose that $L\in U_i\cap \cl_{j'}$ for some $j'$.  Since $U_i\cap \cl_i=\emptyset$, we must have $j'\neq i$.  By (\ref{ueq}), we may write $L=\sum_{j=1}^{d-m}u_j$, with $u_j\in (\cl_j)\cap U_{i-1}$ for $j\neq i$ and $u_i\in (\cl_i)$.  Therefore $L-u_{j'}=\sum_{j\neq j'}u_j$ and 
\begin{equation*}
L-u_{j'}\in (\cl_{j'})\cap \sum_{j\neq j'}(\cl_j).
\end{equation*}
So $L-u_{j'}\in U_0\subset U_{i-1}$.  But $u_{j'}\in U_{i-1}$, which implies that $L\in U_{i-1}$.  This contradicts the assumption that $U_{i-1}\cap \cl=\emptyset$.  Therefore $U_i\cap \cl=\emptyset$.

Let $U_0,\ldots, U_{d-m}$ be vector spaces defined as above.  Let $Y=Z(U_{d-m})$.  We claim that $Y\not\subset |\ch|$, $\cl|_Y$ is a linearly independent set, and $\dim Y=d$.  Since $U_{d-m}\cap \cl=\emptyset$, we have that $Y\not\subset |\ch|$.  As $\dim U_i\cap (\cl_i)=\dim (\cl_i)-1$ for $i>0$, $U_i\subset U_{d-m}$ for all $i$, and $U_{d-m}\cap \cl=\emptyset$, we must have that 
\begin{equation*}
\dim U_{d-m}\cap (\cl_i)=\dim (\cl_i)-1
\end{equation*}
for all $i$.  Therefore 
\begin{equation*}
\dim (U_{d-m}+(\cl_i))=\dim U_{d-m}+\dim(\cl_i)-\dim U_{d-m}\cap (\cl_i)=1+\dim U_{d-m}.
\end{equation*}
Let $\ch_i$ be the set of hyperplanes defined by the elements of $\cl_i$.  The above equation shows that $\ch_i|_Y$ consists of a single hyperplane of $Y$.  So $\ch|_Y$ consists of at most $d-m$ hyperplanes of $Y$.  If $\dim Y=d$, then the fact that $\dim \cap_{H\in H|_Y}H=m$ (since $U_{d-m}\subset (\cl)$), and that $\#\ch|_Y \leq d-m$, imply that $\cl|_Y$ is a linearly independent set.  So it remains to show that $\dim Y=d$, or equivalently, that $\dim U_{d-m}=n-d$.  Repeatedly applying the equation $\dim (U+V)=\dim U+\dim V-\dim U\cap V$ we get that
\begin{equation}
\label{dim1}
\dim \sum_{i=1}^{d-m}(\cl_i)=\dim (\cl)=n-m=\sum_{i=1}^{d-m}\dim (\cl_i)-\sum_{j=1}^{d-m-1}\dim \left((\cl_{j+1})\cap \sum_{i=1}^j(\cl_i)\right)
\end{equation}
and
\begin{multline}
\label{dim2}
\dim U_{d-m}=\dim \sum_{i=1}^{d-m}((\cl_i)\cap U_{d-m})\\
=\sum_{i=1}^{d-m}\dim (\cl_i)\cap U_{d-m}-\sum_{j=1}^{d-m-1}\dim \left((\cl_{j+1})\cap U_{d-m}\cap \sum_{i=1}^jU_{d-m}\cap(\cl_i)\right).
\end{multline}
We claim that
\begin{equation*}
(\cl_{j+1})\cap U_{d-m}\cap \sum_{i=1}^jU_{d-m}\cap (\cl_i)=(\cl_{j+1})\cap \sum_{i=1}^j(\cl_i).
\end{equation*}
One inclusion is trivial.  For the other, let $u\in (\cl_{j+1})\cap \sum_{i=1}^j(\cl_i)$.  This means that we have an equation $u=\sum_{i=1}^ju_i$ where $u\in(\cl_{j+1})$ and $u_i\in (\cl_i)$.  It follows easily from the definition of $U_0$ that $u,u_1,\ldots,u_j\in U_0\subset U_{d-m}$.  The equation then follows.  We also have $\dim (\cl_i)\cap U_{d-m}=\dim (\cl_i)-1$.  Therefore, from Equations (\ref{dim1}) and (\ref{dim2}), we get
\begin{align*}
\dim U_{d-m}&=\sum_{i=1}^{d-m}\dim (\cl_i)-\sum_{j=1}^{d-m-1}\dim \left((\cl_{j+1})\cap \sum_{i=1}^j(\cl_i)\right)
-(d-m)\\
&=n-m-(d-m)=n-d
\end{align*}
as was to be shown.
\end{proof}

\bibliography{integral}

\providecommand{\bysame}{\leavevmode\hbox to3em{\hrulefill}\thinspace}
\providecommand{\MR}{\relax\ifhmode\unskip\space\fi MR }
\providecommand{\MRhref}[2]{%
  \href{http://www.ams.org/mathscinet-getitem?mr=#1}{#2}
}
\providecommand{\href}[2]{#2}
\begin{thebibliography}{1}

\bibitem{Ev}
J.-H. Evertse and K.~Gy{\H{o}}ry, \emph{Finiteness criteria for decomposable
  form equations}, Acta Arith. \textbf{50} (1988), no.~4, 357--379.

\bibitem{Fu2}
Hirotaka Fujimoto, \emph{Extensions of the big {P}icard's theorem}, T\^ohoku
  Math. J. (2) \textbf{24} (1972), 415--422.

\bibitem{Gr}
Mark~L. Green, \emph{Holomorphic maps into complex projective space omitting
  hyperplanes}, Trans. Amer. Math. Soc. \textbf{169} (1972), 89--103.

\bibitem{No}
Junjiro Noguchi and J{\"o}rg Winkelmann, \emph{Holomorphic curves and integral
  points off divisors}, Math. Z. \textbf{239} (2002), no.~3, 593--610.

\bibitem{Ru2}
Min Ru, \emph{Geometric and arithmetic aspects of {$\bold P\sp n$} minus
  hyperplanes}, Amer. J. Math. \textbf{117} (1995), no.~2, 307--321.

\end{thebibliography}
\end{document}